\documentclass[11pt]{amsart}

\usepackage{amssymb, array, verbatim, amscd}
\usepackage{amsmath, amsfonts}

\theoremstyle{plain}
\newtheorem{theorem}{Theorem}[section]
\newtheorem{corollary}[theorem]{Corollary}
\newtheorem{lemma}[theorem]{Lemma}
\newtheorem{proposition}[theorem]{Proposition}

\newcounter{abc}

\newcounter{ABC}

\newcommand{\Z}{{\mathbb Z}}

\newcommand{\E}{{\mathcal E}}
\newcommand{\A}{{\mathcal A}}

\begin{document}
\title[Splitting the Automorphism Group]{Splitting the Automorphism Group of an Abelian p-Group }
\author{Maria Alicia Avi\~no}
\address[Avi\~no]{Department of Mathematical Sciences,\\ New Mexico State University,\\
Las Cruces, 88003, }
\email{mavino@nmsu.edu}

\begin{abstract}
Let G be an abelian p-group sum of finite homocyclic groups $G_i$. Here, we
determine in which cases the automorphism group of G splits over ker$\sigma$, where $\sigma : $Aut$(G) \rightarrow \prod _i $Aut$(G_i/pG_i)$ is the natural epimorphism.
\end{abstract}
\maketitle
\section{Preliminaries}
Throghout this paper,  $p$ is an arbitrary prime  and $r$ is a fixed ordinal number.
 Let $G=\oplus_{i\leq r} G_i$ be an abelian p-group, such that
$G_i$ is an homocyclic group of exponent $p^{n_i}$ and finite p-rank $r_i$, with
$n_i<n_{i+1}$ for all $i$.
It is known that $ \E$ the endomorphism ring of G is isomorphic
to the ring $E(G)$ of all row finite $r\times r$-matrices
$(A_{ij})$ where $A_{ij}\in Hom (G_i,G_j)$.
We denote by ${\A}(H)$, the automorphism group of a group $H$ and consider
${\A}(H)$ as the group of units of the endomorphism ring $E(H)$.

Let $\sigma $ be the natural epimorphism of ${\A}( G)$ onto the product
of the ${\A}(G_i/pG_i)$. We have the following exact sequence (see \cite {2}, page 256)
$$1\rightarrow ker \sigma \rightarrow {\A}( G) \rightarrow
\prod _i{\A}(G_i/pG_i)\cong \prod _i  GL_{r_i}({\bf Z}_p)\rightarrow 1 ,\leqno (1)$$
where  $GL_{r_i}({\bf Z}_p)$ is the general linear
group of $r_i\times r_i$-matrices over the field ${\bf Z}_p$.

In this paper we prove a Theorem which together with Theorem 1.1
proved in \cite{1} and \cite{3}, give a necessary and sufficient condition for the decomposition
 of ${ \A}(G)$ as a semidirect product of ker$\sigma $=$\Delta (G)$ by
$\Pi (G)=\prod _i  GL_{r_i}({\bf Z}_p)$, whenever
$p\ge 5$. For the cases $p=2,3$ we
give sufficient and necessary conditions for such decomposition in case
$n_i+1<n_{i+1}$ for all $i$.

 Because the p-rank of $G_i$
 is finite,  $Hom (G_i,G_j)\cong  p^{n_j-n_i}M_{r_i\times r_j}
({\bf Z}_{p^{n_i}})$ for $i<j$, where $M_{r_i\times r_j}
({\bf Z}_{p^{n_i}})$ denotes the additive group of $r_i\times r_j-$matrices over the
integers modulo $p^{n_i}$.

Evidently  $\Delta (G)=1+I$ where $$I=\{(A_{ij})_{r\times r}\in E(G)\vert A_{ii}\equiv 0 \hbox
{(mod $p$),  for all } i \}.$$
If $r$ is a natural number, then $G$ is finite and because
$GL_{r_{i}}({\bf Z}_{p})$ does not have any normal p-subgroup,
$\Delta (G)$ is the maximal normal p-subgroup of ${\A}(G)$, denoted $O_p({\A}(G))$.

Consider  the exact sequences
$$1 \rightarrow ker \lambda_i \rightarrow {\A}( G_i) \rightarrow
GL_{r_{i}}({\bf Z}_{p}) \rightarrow 1\leqno {(2i)}$$
where $GL_{r_{i}}({\bf Z}_{p})$ is the general linear group of $r_{i} \times r_{i}$-matrices
over the field ${\bf Z}_{p}$.
It follows that $ker\lambda_i=J($End $G_i) +1_i=O_p({\A} G_i)=\Delta (G_i)$.

By Theorem I, II.2 and III.4 in [1] and the Theorem in [3] we have the following.

\begin{theorem}

 1.1.1) For $p\ge 5$, the exact sequence
(2i) splits if and only if $G_{i}$ is of type either $(p^{n_{i}})$
or $(p,...,p)$ with $n_{i}>1$.

1.1.2) For $p=3$, the exact sequence (2i)
splits if and only if $G_{i}$ is of type  $(3^{n_i})$,
$(3^{n_i},3^{n_i})$ or $(3,...,3)$ with $n_i>1$.

1.1.3) For $p=2$ the exact sequence (2i)
splits if and only if $G_{i}$ is of type  $(2^{n_i})$,
$(2^{n_i},2^{n_i})$, $(2^{n_i},2^{n_i},2^{n_i})$ or $(2,...,2)$ with $n_i>1$.

\end{theorem}
\section{ ${\A}(G)$ as semidirect product of $\Delta (G)$ and $\Pi (G)$ }

This section is devoted to the proof of Theorem 2.1 .

\begin{theorem}

 2.1.1) For $p\ge 5$, the exact sequence $(1)$ splits if and only if
 the exact sequences $(2i)$ split for all $i$.

2.1.2) If $n_{i+1}-n_i>1$ for all $i$, the exact sequence $(1)$
splits if and only if the exact sequences $(2i)$ split for all $i$.
\end{theorem}
This Theorem is a consequence of the following lemmas and propositions.

\begin{proposition} If the exact sequences $(2i)$ split for
 all $i$, then the exact sequence $(1)$ splits.
\end{proposition}

\begin{proof}
 Because the exact sequences $(2_i)$ split for all $i$,
there exists a family of morphisms $\{\theta_{i}: GL_{r_{i}}({\bf Z}_{p})
\rightarrow {\A} (G_{i})\}_{r}$ such that $\theta_{i} \lambda_{i}=1_i$,
where $1_i$ is the identity of ${\A}( G_i)$. So we have the diagram
\[
\begin{array}{ccc}
\A( G ) &{\overset{\sigma}{\longrightarrow}}  & \prod_i GL_{r_i}(\Z_p)\\
{\iota}{\uparrow } & {}&{\vert}{\uparrow}  \\
\prod_i {\A} (G_i)& {\overset{\sigma \vert _\Pi}{\longrightarrow}}&
\prod_{i} GL_{r_{i}}(\Z_{p})\\
 {p_i}{\downarrow }& {} &{\delta _i} {\downarrow }\\
 {\A} (G_{i}) &{}^{\overset{\lambda_{i}}{\longrightarrow}}_{
\underset{\theta_{i}}{\longleftarrow}} &
 GL_{r_i}(\Z_{p}) \end{array}
\]

where $\iota  $ is the canonical monomorphism, $p_i$, $\delta _i$
are the canonical epimorphisms and $\sigma \vert _\Pi$ is the restriction of $\sigma $
to $\prod_r {\A} (G_i)$.
Now consider the morphisms $f_{i}=\delta_{i}\theta _i$, for all
$i $.

There exists $\theta: \prod_{i}GL_{r_{i}}(\Z_{p})
 \rightarrow \prod _{r}{\A} (G_i)$
such that $\theta p_i=f_{i} $ for all $i$,
 and the diagram commutes.
It is clear that $\theta \sigma\vert_{\Pi}  =1_{\Pi}$ where $1_{\Pi}$
 is the identity map of $\prod_{i} GL_{r_{i}}(\Z_{p})$,
then $\theta \sigma = 1_{\Pi}$ so the Proposition  holds.
\end{proof}

In the following lemma we consider the  exact sequence $(1)$ for the
  group $pG$.
$$
1 \rightarrow \Delta (pG) \rightarrow { \A}(pG) {\buildrel\sigma'
\over \rightarrow}\Pi (pG) \rightarrow 1
\leqno {(3)}
$$
\begin{lemma} If the exact sequence $(1)$ splits then the
exact sequence $(3)$ splits.
\end{lemma}

\begin{proof} Set $\alpha : G \rightarrow pG$ such that $g \alpha =pg.$
But  $ker \alpha = G[p]$ is the socle of G.
 This morphism induces a morphism $\tau : {\A}(G) \rightarrow {\A} (pG)$.
 Because $\Delta (G)$ is a characteristic subgroup $\Delta (G)\tau \subset \Delta (pG)$ so we have the following exact commutative diagram

\[
\begin{array}{ccccccccc}
1&\rightarrow & \Delta (G) & \rightarrow & {\A} (G) &{\overset {\sigma}
{\rightarrow}} & \Pi (G) & \rightarrow & 1\\
&{ } & \downarrow & { } &{^{\tau}{\downarrow}} & {}
&{^{\underline\tau }{\downarrow}}& {}&{}\\
1&\rightarrow & \Delta (pG) & \rightarrow & {\A}(pG) &
\rightarrow & \Pi (pG) & \rightarrow & 1\end{array}
\]
and two cases

a) If $n_1>1$, then $G[p]\subset pG $, so $\Pi (G)=\Pi (pG)$
and $\underline \tau$  is the identity. Thus if $\theta$ splits
$(1)$ then $\theta \tau $ splits $(3)$.

b) If $n_1=1$, then $G_1$ is of type $(p,p,...,p)$ and
$G=G_1\oplus \overline G_2$ where
$\overline G_2=\oplus_{i\ge 2} G_i$.
Therefore $\Pi (pG)=\Pi (\overline G_2)$.
Thus, if $\iota $ is the canonical monomorphism of $\Pi (pG)$ to
$\Pi (G) $ we have that if $\theta $ splits $(1)$, then the morphism
$\iota \theta \tau   $ splits $(3)$ and the Lemma holds.
\end{proof}

\begin{lemma}  Set $\overline G_2=\oplus_{i\ge 2} G_i$. If $n_1=1$,
then $(1)$ splits if and only if the exact sequence
$1\rightarrow \Delta (\overline G_2) \rightarrow {\A} (\overline G_2)
{\buildrel \sigma_2 \over \rightarrow} \Pi(\overline G_2)\rightarrow 1$
splits.
\end{lemma}

\begin{proof} Let $g\in G$ and $g=g_1+g_2$, $g_1\in G_1$, $g_2\in \overline G_2$.
Then $\gamma :G\rightarrow \overline G_2$ such that $ g\gamma=g_2$
is a morphism with $ker \gamma =G_1\subseteq G[p]$ the socle of $G$.
This morphism induces a morphism
$\overline  \gamma : {\A} (G) \rightarrow \A (\overline G_2)$.

Here $\iota _2$ is the natural inclusion and
$\overline p(A_{ii})_r=(A_{ii})_{i>1}$. Then the following
 diagram commutes
\[
\begin{array}{ccc}
{\A} (G )&{\overset{\sigma}{\longrightarrow}} & \Pi(G)
\\{^{\iota_2}
{\uparrow}}{^{\downarrow}{\overline\gamma}}
& &{^{\iota}{\uparrow}}{^{\downarrow}{\overline p}}\\
 {\A}(\overline G_2) & {\overset {\sigma_2}{\rightarrow}} & \Pi (\overline G_2) \end{array}
\]
and our claim holds.
\end{proof}
\begin{lemma}
 Let $n_1=2$ and $n_2\ge 4$. If $(1)$ splits then
$(2_1)$ splits.
\end{lemma}
\begin{proof}
 Set $A=(A_{ij})\in {\A}(G)$, we can define a map
$\mu : {\A}( G) \rightarrow {\A} (G_1)$ so that $\mu (A)=A_{11}$.
In this case, the map $\mu $ is a morphism. In fact, if $B=(B_{ij})\in {\A}(G)$ then  $(BA)_{11}=B_{11}A_{11}$ since $A_{i1}\equiv 0 \hbox { (mod } p^{n_i-n_1})$
and $n_i-n_1\ge 2$ for $i>1$.

Now, consider the diagram
\[
\begin{array}{ccc}
{\A} (G) & {\overset{\underline\sigma}{\rightarrow}} & \Pi(G) \\
\downarrow{\mu } &&{\iota_1}{\uparrow }{\downarrow}{\mu'}\cr
 {\A} (G_1) &{\overset {\lambda_1}{\rightarrow}} & \Pi (G_1)
\end{array}
\]
where $\mu'$ is the morphism induced by $\mu$, and $\iota_1$
is the canonical monomorphism. If $\theta$ splits $(1)$
then $\iota_1 \theta \mu$ splits $(2_1)$.
\end{proof}

\begin{proposition} If  $n_{i+1}-n_i>1$
for all $i$ and the exact sequence $(1)$ splits then  $(2_i)$ split
for all $i$.
\end{proposition}

\begin{proof}

 Suppose there exists $\theta$ such that $(1)$ splits,
but there exists  $i$ such that $(2_i)$ does not split.
Let $t$  be the least $i$ such that  $(2_i)$ does not split.
We have two cases

a) $t>1$, so $n_t> 2$.
Applying the Lemma 2.3 we can obtain that the exact sequence
$$1\rightarrow \Delta ({\underline G}) \rightarrow \A({\underline G})
 \rightarrow \Pi ({\underline G}) \rightarrow 1 $$
splits where ${\underline G}= p^{n_t-2}G=\oplus_{i\ge t}\underline G_i$, with
$\underline G_i=p^{n_t-2} G_i$ for
$i\ge t$.
But the exponent of ${\underline G_t}$ is $p^{2}$ and the exponent of
$G_{t+1}$ is $\ge p^4$.

Then using Lemma 2.5 we obtain that the exact
sequence
$$1\rightarrow \Delta (\underline G_t) \rightarrow {\A}(\underline G_t)
\rightarrow \Pi (\underline G_t) \rightarrow 1 $$
splits. But this is impossible because the rank
of $\underline G_t$ is the same as that  of $G_t$,
that is  $r_t\ge 2$ if $p\ge 5$, $r_t\ge 3$
if $p=3$  and  $r_{t}\ge 4$ if $p=2$, contradicting Theorem 1.1.

b) If $t=1$, but $n_1\ge 2$, we proceed as in a).
If $t=1$ but $n_1=1$,  $G_1$
is an elementary abelian group and by Theorem 1.1, $(2_1)$ splits,
 which is not the case.
\end{proof}

\begin{proposition} If $p\ge 5$ and  the exact sequence $(1)$ splits then  the
exact sequences $(2_i)$
split, for all $ i$.
\end{proposition}

\begin{proof}
Suppose that there exists t such that $(2_t)$ does not split. Then by Theorem 1.1, $G_{t}$ is a subgroup of exponent $p^{n_t}$,
$n_t>1$ and p-rank $r_t>1$.
Using Lemmas 2.3 and 2.4, we may assume  $t=1$ and $n_t=2$. Thus
we may assume that $(1)$ splits and $(2_1)$ does not split with
$n_1=2$ and  $r_1\ge 2$.

Let $\overline A=1_\Pi +(E_{1r_1},0,...,0)\in \Pi (G)$, where
$E_{1r_1}=(e_{uv})\in GL_{r_1}(\Z_p)$ such that $e_{1r_1}=1$ and
 $e_{uv}=0$ for $(u,v)\ne (1,r_1)$, and $A=1+B=1+(B_{ij})\in \A(G)$,
 where $B_{11}=E_{1r-1}$ and $B_{ij}=0$ for $(i,j)\ne (1,1)$.

If $\theta $ is the monomorphism such that $\theta \sigma = 1_\Pi $,
 then $\theta (\overline A)\equiv A$ (mod $\Delta (G)$). So
$\theta (\overline A)=A(1+\overline C)$ with
$1+\overline C \in \Delta (G)=1+I$. Set
$A\overline C=C=(c_{uv})\in I$, then we have that $ \overline A \theta=A+C$.

Because the order of $\overline{A}$ is $p$,  $1=(A+C)^p$.
Therefore
 $$1=[1+(B+C)]^p= 1+p(B+C)+{p\choose 2 }(B+C)^2+{p\choose 3}(B+C)^3+$$
$$.....+{p\choose p-1}(B+C)^{p-1}+(B+C)^p$$

Since $(B+C)^2\in I$, we have that the block $(1,1)$ of
$(B+C)^t$ is $0$ for $t\ge 4$ and for the matrices $(B+C)^2$, $(B+C)^3$
 the block $(1,1)$ is congruent to $p$ modulo $p^2$.
For $p\ge 5$, we have in the place $(1,r_1)$  \[0\equiv
p(1+pc_{1r_1})\hbox { (mod }p^2 ),\] which is impossible and the
proposition holds.
\end{proof}
\begin{corollary}

 2.8.1) For $p\ge 5$, ${\A} (G) $ is a semidirect product of $\Delta (G)$ by
$\Pi (G)$ if and only if ${\A}(G_{i})$ is a semidirect product
of $O_{p}({\A} (G_{i}))$ by $GL_{r_{i}}(\Z_{p})$.

2.8.2) For $p=3$, $p=2$ and $n_{i+1}-n_{i}>1$, $i$
${\A} (G)$ is a semidirect product of $\Delta (G)$ by
$\Pi (G)$ if and only if ${\A} (G_{i})$ is a semidirect product of
$O_{p}({\A} (G_{i}))$ by $GL_{r_{i}}(\Z_{p})$.
\end{corollary}

\begin{corollary}

 ${\A} (G)$ is a semidirect product of $\Delta (G)$ by
$\Pi (G)$ if and only if

2.9.1) For $p\ge 5$,  $G_i$ has rank $1$ for $i>1$ and  $G_1$ either
is elementary abelian or  has  rank $1$ if $n_1>1$.

2.9.2) For $p=3$ and $n_{i+1}-n_i>1$, for all $i$, $G_i$ has rank $\le 2$ for $i>1$
and  $G_1$ either is elementary abelian or  has  rank $\le 2$ if $n_1>1$.

2.9.3) For $p=2$ and $n_{i+1}-n_{i}>1$ for all $i$, $G_i$ has rank
$\le 3$ for $i>1$ and  $G_1$ either is elementary abelian or  has
  rank $\le 3$ if $n_1>1$.
\end{corollary}


\begin{thebibliography}{99}



\bibitem [1]{1} Beletskii, P and Morgado, E. {\it On the automorphism group of a finite abelian
 p-group.} Ukrain. Mat. Zh. {\bf 33} No 5, 589-596.

\bibitem [2]{2}  Fuchs, L. {\it Infinite Abelian Groups}.
 Academic Press {\bf Vol II} (1973).

\bibitem [3]{3}  Morgado, E. {\it On the automorphism group of a finite abelian
 p-group.} Ukrain. Mat. Zh. {\bf 32} No 5, 617-622.

\end{thebibliography}
\end{document}